\documentstyle[12pt, fullpage]{amsart}
\begin{document}
\newcommand{\createtitle}[2]{\title{#1}\author{Greg Martin}\address{Department of Mathematics\\University of Toronto\\Canada M5S 3G3}\email{gerg@@math.toronto.edu}\subjclass{#2}\maketitle}
\newcommand{\make}[1]{\label{#1sec}\noindent\input{#1}}
\newcommand{\Mmake}[1]{\label{#1sec}\noindent}


\newtheorem{theorem}{Theorem}
\newtheorem{lemma}[theorem]{Lemma}
\newtheorem{corollary}{Corollary}[theorem]
\newtheorem{proposition}[theorem]{Proposition}

\newenvironment{pflike}[1]{\noindent{\bf #1}}{\vskip10pt} 
\newenvironment{proof}{\begin{pflike}{Proof:}}{\qed\end{pflike}}

\def\({\left(}\def\){\right)}

\newcommand{\2}[1]{\ifmmode{\cal#1}\else$\cal#1$\fi}
\newcommand{\3}[1]{\#\{#1\}}
\newcommand{\abs}[1]{\left|#1\right|} 
\newcommand{\floor}[1]{\lfloor#1\rfloor}
\newcommand{\bfloor}[1]{\big\lfloor#1\big\rfloor}
\newcommand{\bbfloor}[1]{\bigg\lfloor#1\bigg\rfloor}
\newcommand{\ceil}[1]{\lceil#1\rceil}
\newcommand{\bceil}[1]{\big\lceil#1\big\rceil}
\newcommand{\bbceil}[1]{\bigg\lceil#1\bigg\rceil}

\renewcommand{\mod}[1]{{\ifmmode\text{\rm\ (mod~$#1$)}
 \else\discretionary{}{}{\hbox{ }}\rm(mod~$#1$)\fi}}
\newcommand{\ep}{\varepsilon}
\renewcommand{\implies}{\Rightarrow}
\newcommand{\rmif}{{\rm if\ }}

\newcommand{\half}{{\mathchoice{\textstyle\frac12}{1/2}{1/2}{1/2}}}
\newsymbol\dnd 232D 
\newcommand{\exdiv}{\mathrel{\mid\mid}}

\renewcommand{\lg}[1]{\mathop{\log_{#1}}}
\def\lgs#1^#2{\mathop{\log_{#1}^{#2}}}
\newcommand{\li}{\mathop{\rm li}}

\newcommand{\doublespace}{
  \baselineskip=24pt}
\newcommand{\spaceandahalf}{\parskip6pt\baselineskip=18pt}

\newcommand{\scroll}[1]{\scrollmode#1\errorstopmode}
\newcommand{\comment}[1]{}

\vfuzz=2pt 
\createtitle{Uniform Bounds for the Least Almost-Prime Primitive~Root}{11N69}

\def\implies{\Rightarrow}
\def\lr{$\lambda$\discretionary{-}{}{-}root}
\def\nup{{m(p)}}
\def\Zx#1{{\bf Z}_{#1}^\times}
\def\Zqx{\Zx q}
\def\ord#1#2{{\rm ord}_#1(#2)}
\def\A{{\cal A}}\def\B{{\cal B}}\def\C{{\cal C}}
\def\pth{$p^{\rm th}$}

\section{Introduction}\noindent
A recurring theme in number theory is that multiplicative and additive
properties of integers are more or less independent of each other, the
classical result in this vein being Dirichlet's theorem on primes in
arithmetic progressions. Since the set of primitive roots to a given
modulus is a union of arithmetic progressions, it is natural to study
the distribution of prime primitive roots. Results concerning upper
bounds for the least prime primitive root to a given modulus $q$,
which we denote by $g^*(q)$, have hitherto been of three types. There
are conditional bounds: assuming the Generalized Riemann
Hypothesis, Shoup \cite{Sho:SfPRiFF} has shown that
\begin{equation*}
g^*(q) \ll \big(\omega(\phi(q))\log2\omega(\phi(q))\big)^4(\log q)^2,
\end{equation*}
where $\omega(n)$ denotes the number of distinct prime factors of
$n$. There are also upper bounds that hold for almost all moduli
$q$. For instance, one can
show \cite{Mar:TLPPRatSS} that for all but $O(Y^\ep)$ primes up to
$Y$, we have
\begin{equation*}
g^*(p) \ll (\log p)^{C(\ep)}
\end{equation*}
for some positive constant $C(\ep)$. Finally, one can unimaginatively
apply a uniform upper bound for the least prime in a single arithmetic
progression. The best uniform result of this type, due to Heath-Brown
\cite{Hea:Z-fRfDL-fetc}, implies that $g^*(q) \ll q^{5.5}$.  However,
there is not at present any stronger unconditional upper bound for
$g^*(q)$ that holds uniformly for all moduli $q$. The purpose of this
paper is to provide such an upper bound, at least for primitive roots
that are ``almost prime''.

The methods herein will actually apply for any modulus $q$, not just
those $q$ whose group $\Zqx$ of reduced residue classes is cyclic
(which occurs exactly when $q=2$, 4, an odd prime power, or twice an
odd prime power). We say that an integer $n$, coprime to $q$, is a
{\it\lr}\mod q if it has maximal order in $\Zqx$. We see that \lr{}s
are generalizations of primitive roots, and we extend the notation
$g^*(q)$ to represent the least prime \lr\mod q for any integer
$q\ge2$.

We also recall that a $P_k$ integer is one that has at most $k$ prime
factors, counted with multiplicity. For any integer $k\ge 1$, we let
$g^*_k(q)$ denote the least $P_k$ \lr\mod q (so that
$g^*_1(q)=g^*(q)$). We may now state our main theorem.

\begin{theorem}For all integers $q$, $r\ge2$ and all $\ep>0$, we have
\begin{equation*}
g^*_r(q) \ll_\ep q_c^{1/4+1/(4(r-1-\delta_r))+\ep},
\end{equation*}
where $q_c$ is the largest odd cubefree divisor of $q$, and one can
take
\begin{equation*}
\delta_2=0.0044560, \quad \delta_3=0.074267, \quad \delta_4=0.103974,
\end{equation*}
and $\delta_r=0.1249$ for any $r\ge5$.
%
%
\label{mainthm}
\end{theorem}

In the context of primitive roots, Theorem \ref{mainthm} can be
improved to the extent of replacing the largest cubefree divisor of
the modulus with the largest squarefree divisor. The result is as
follows:

\begin{theorem}
Let $p$ be an odd prime, and let $q$ be a power of $p$ or twice a
power of $p$. Then
\begin{equation*}
\begin{split}
g^*_2(q) &\ll p^{1/2+1/873}; \\
g^*_3(q) &\ll p^{3/8+1/207}; \\
g^*_4(q) &\ll p^{1/3+1/334}; \\
g^*_r(q) &\ll_r p^{1/4+O(1/r)}.
\end{split}
\end{equation*}
\label{primrootthm}
\end{theorem}

\noindent (The exponents here are simply approximations to the
corresponding exponents in Theorem \ref{mainthm}.) By comparison, from
the work of Mikawa on small $P_2$s in almost all arithmetic
progressions \cite{Mik:A-PiAPaSI}, one can easily derive that
\begin{equation*}
g^*_2(q)\ll q(\log q)^5{\phi(q)\over\phi(\phi(q))},
\end{equation*}
which is majorized by the above theorems. We remark that the \lr{s} we
find to establish Theorems \ref{mainthm} and \ref{primrootthm} are
squarefree and have no small prime factors (``small'' here meaning up
to some fixed power of $q_c$).

The proof of Theorem \ref{mainthm} uses the weighted linear sieve,
specifically results due to Greaves \cite{Gre:AWSoBT},
\cite{Gre:TWLSaSetc.} (see equation (\ref{Gresult}) below). We note
that conjecturally, there is some choice of weight function $W$ in the
weighted linear sieve which would allow us to take $\delta_r$
arbitrarily small in Theorem \ref{mainthm}; this would also allow us
to replace the first three exponents in Theorem \ref{primrootthm} by
$1/2+\ep$, $3/8+\ep$, and $1/3+\ep$ respectively. We also note that if
the generalized Lindel\"of hypothesis for the $L$-functions
corresponding to certain characters (the ones in the subgroup $G$
defined in Lemma \ref{charlincomb} below) were true, we could employ
much stronger character sum estimates than Lemma \ref{Cqbd} below,
allowing us to improve Theorem \ref{mainthm} to $g_2^*(q)\ll_\ep
q_c^\ep$.

We would of course like to be able to show the existence of small
prime primitive roots rather than $P_2$ primitive roots. In his work
on the analogous problem of finding primes in arithmetic progressions,
Heath-Brown \cite{Hea:SZatLPiaAP} first treats the case where the
$L$-function corresponding to a real Dirichlet character has a real
zero very close to $s=1$. Although it is certainly believed that these
``Siegel zeros'' do not exist, disposing of this case allowed
Heath-Brown in \cite{Hea:Z-fRfDL-fetc} to work with a better zero-free
region for Dirichlet $L$-functions than is known unconditionally.

Similarly, if we assume the existence of a sufficiently extreme Siegel
zero, we can show the existence of small prime primitive roots, as the
following theorem asserts.

\begin{theorem}Let $\ep>0$, let $q$ be an odd prime power or
twice an odd prime power, and let $\chi_1$ denote the nonprincipal
quadratic Dirichlet character\mod q. Suppose that $L(s,\chi_1)$ has a
real zero $\beta$ satisfying
\begin{equation*}
\beta > 1-\frac1{A(\ep)\log q},
\end{equation*}
where $A(\ep)$ is some sufficiently large real number depending on
$\ep$ (but not on $q$). Then
\begin{equation*}
g^*(q) \ll_\ep p^{3/4+\ep}.
\end{equation*}
\label{Siegelthm}
\end{theorem}

It is plausible that in the weighted linear sieve, one can derive a
corresponding upper bound of the same order of magnitude as the lower
bound. In this case, the exponent $3/4+\ep$ in Theorem \ref{Siegelthm}
could be replaced by $1/2+1/873$, the exponent associated with
$g_2^*(q)$ in Theorem~\ref{primrootthm}.

It is a pleasure to thank Trevor Wooley and Hugh Montgomery for their
many helpful suggestions on improving this paper, and for their
guidance in general. The author would also like to thank Andrew
Granville, George Greaves, Ram Murty, and Amora Nongkynrih for
valuable comments regarding recent progress in areas relevant to this
paper. This material is based upon work supported under a National
Science Foundation Graduate Research Fellowship. 

For the remainder of this paper, the constants implicit in the $\ll$
and $O$-notations may depend on $\ep$, $\eta$, and $r$ where
appropriate. We denote the cardinality of the set $S$ by $\abs S$.  As
before, $\omega(n)$ is the number of distinct prime factors of $n$,
and $\Omega(n)$ is the number of prime factors of $n$ counted with
multiplicity, so that an integer $n$ is a $P_r$ precisely when
$\Omega(n)\le r$.

\section{The characteristic function of \lr{}s}\noindent
Let $\gamma$ be the characteristic function of those integers that are
{\lr}s\mod q. Since $\gamma$ is periodic with period $q$ and is
supported on reduced residue classes, $\gamma$ can be written as a
linear combination of Dirichlet characters. The following lemma
exhibits this linear combination explicitly. Let $E(q)$ denote the
exponent of the group $\Zqx$, so that an integer $n$, coprime to $q$,
is a \lr\mod q precisely when the multiplicative order of $n$ is
$E(q)$. Let $S(q)$ denote the largest squarefree divisor of $E(q)$. We
notice that $\phi(q)$, $E(q)$, and $S(q)$ all have exactly the same
prime divisors.

\begin{lemma}Let $G$ be the subgroup of characters\mod q given by
\begin{equation*}
G = \{\chi^{E(q)/S(q)}: \chi\mod q\}.
\end{equation*}
For every prime $p$ dividing $\phi(q)$, let $\nup$ denote the number
of independent characters of order $p$ in $G$. For every character
$\chi\mod q$, let $\sigma(\chi)$ denote the order of $\chi$, and
define
\begin{equation}
c_\chi = 
\begin{cases}{\displaystyle\prod_{p\mid\sigma(\chi)}
\(\frac{-1}{p^{\nup}}\) \prod\begin{Sb}p\mid\phi(q) \\
p\dnd\sigma(\chi)\end{Sb} \( 1-\frac1{p^{\nup}} \)} &\rmif\chi\in G,\\
0 &\text{\rm otherwise.}\end{cases}
\label{cchi}
\end{equation}
Then, for any integer $n$,
\begin{equation}
\gamma(n) = \sum_{\chi\mod q} c_\chi\chi(n). \label{gamma}
\end{equation}
\label{charlincomb}
\end{lemma}

\begin{pflike}{Remarks:}
For simplicity we write $c_0$ for $c_{\chi_0}$. We note that the
exponent of $G$ is $S(q)$, the number of characters in $G$ whose order
divides $p$ is $p^{\nup}$, and the number of characters in $G$ whose
order equals $d$ is
\begin{equation}
\prod_{p\mid d} \( p^\nup-1 \).
\label{orderisd}
\end{equation}
We note that $c_0$ is the probability that a randomly chosen element
of $\Zqx$ has order $E(q)$, i.e., the number of \lr{}s\mod q less than
$q$ is $c_0\phi(q)$. We also note that $\Zqx$ is cyclic if and only if
$E(q)=\phi(q)$. When $\Zqx$ is cyclic, the definition (\ref{cchi})
of $c_\chi$ reduces to
\begin{equation*}
c_\chi = {\phi(\phi(q))\over\phi(q)}
{\mu(\sigma(\chi))\over\phi(\sigma(\chi))}.
\label{cchicyc}
\end{equation*}
In particular, $c_0=\phi(\phi(q))/\phi(q)$ in this case.
\end{pflike}

\begin{proof} 
The lemma clearly holds for $(n,q)>1$, since both sides of equation
(\ref{gamma}) are zero, and thus for the remainder of the proof, we
assume that $(n,q)=1$. From the standard properties of characters,
for every prime $p$ dividing $\phi(q)$ we have
\begin{equation*}
\sum\begin{Sb}\chi\in G \\ \chi^p=\chi_0\end{Sb}
\chi(n)=\begin{cases}p^{\nup} &\rmif n^{E(q)/p}\equiv1\mod q,
\\0 & \text{otherwise,}\end{cases}
\end{equation*}
since the number of characters being summed over is $p^{\nup}$, as
noted above. We rewrite this as
\begin{equation}
\( 1 - \frac1{p^\nup}\)\chi_0(n) - \frac1{p^\nup}\sum\begin{Sb}\chi\in
G \\ \chi^p=\chi_0 \\
\chi\ne\chi_0\end{Sb}\chi(n)=\begin{cases}1&\rmif n^{E(q)/p}\not\equiv1\mod
q,\\0&\text{otherwise.}\end{cases}
\label{psum2}
\end{equation}
Now $n$ is a \lr{} if and only if for every prime $p$ dividing
$\phi(q)$, we have $n^{E(q)/p}\not\equiv1\mod q$. Therefore, by using the
relation (\ref{psum2}) for all primes dividing $\phi(q)$, we see that
\begin{equation*}
\prod_{p\mid\phi(q)} \Bigg( \( 1 - \frac1{p^\nup}\)\chi_0(n) - \frac1{p^\nup}
\sum\begin{Sb}
\chi\in G \\ \chi^p=\chi_0 \\ \chi\ne\chi_0
\end{Sb}\chi(n) \Bigg) = \gamma(n).
\end{equation*}
When we expand this product, the only characters that appear are
those in $G$. For such a character $\chi$, the coefficient is
\begin{equation*}
\prod_{p\mid\sigma(\chi)} \( \frac{-1}{p^{\nup}} \)
\prod\begin{Sb}
p\mid\phi(q) \\ p\dnd\sigma(\chi)
\end{Sb}\( 1-\frac1{p^{\nup}} \),
\end{equation*}
which is the same as the definition (\ref{cchi}) of $c_\chi$. This
establishes the lemma.
\end{proof}

\begin{lemma}Let $c_\chi$ be defined as in (\ref{cchi}). Then
\begin{equation}
\sum_{\chi\mod q} \abs{c_\chi} = 2^{\omega(\phi(q))}c_0.
\label{firstass}
\end{equation}
In particular, for any $\ep>0$,
\begin{equation}
\sum_{\chi\mod q} \abs{c_\chi} \ll (S(q))^\ep.
\label{cchisumbd}
\end{equation}
\label{cchisum}
\end{lemma}

\begin{proof}
From the definition (\ref{cchi}) of $c_\chi$, we have
\begin{equation*}
\sum_{\chi\mod q} \abs{c_\chi} = \sum_{d\mid S(q)} \Bigg|\prod_{p\mid
d} \(\frac{-1}{p^{\nup}}\) \prod\begin{Sb}p\mid\phi(q) \\ p\dnd
d\end{Sb} \( 1-\frac1{p^{\nup}} \)\Bigg| \sum \begin{Sb}\chi\in G \\
\sigma(\chi)=d\end{Sb} 1.
\end{equation*}
Using the remark (\ref{orderisd}) to evaluate the inner sum, this
becomes
\begin{equation*}
\begin{split}
\sum_{\chi\mod q} \abs{c_\chi} &= \sum_{d\mid S(q)} \prod_{p\mid
d} \(\frac1{p^{\nup}}\) \prod\begin{Sb}p\mid\phi(q) \\ p\dnd
d\end{Sb} \( 1-\frac1{p^{\nup}} \) \prod_{p\mid d} \(
p^\nup-1 \) \\
&= \prod_{p\mid\phi(q)} \( 1-\frac1{p^{\nup}} \) \sum_{d\mid S(q)} 1 \\
&= c_0\cdot2^{\omega(S(q))}.
\end{split}
\end{equation*}
This establishes the first assertion (\ref{firstass}) of the lemma,
since $\omega(S(q))=\omega(\phi(q))$, and also the second assertion
(\ref{cchisumbd}), since $c_0\le1$ and $2^{\omega(n)}\ll n^\ep$ for
all $n$.
\end{proof}

Let $q_c$ denote the largest odd cubefree divisor of $q$. For any
prime $p$ and any nonzero integer $n$, let $\ord pn$ denote the
largest integer $r$ such that $p^r$ divides $n$. Let
$\alpha=\max\{3,\ord2q\}$, and let $\tilde{q_c}=2^\alpha q_c$. Thus
$\tilde{q_c}$ is almost the largest cubefree divisor of $q$, except
that we allow 8 to divide $\tilde{q_c}$ if it divides $q$.

\begin{lemma}
Every $\chi\in G$ is induced by a character\mod{\tilde{q_c}}.
\label{induced}
\end{lemma}

\begin{proof}
Since every $\chi\in G$ is the $(E(q)/S(q))$-th power of a
character\mod q, it suffices to show the following: for every
character $\chi\mod q$, $\chi^{E(q)/S(q)}$ is periodic with period
dividing $\tilde{q_c}$. For the remainder of the proof, let $\chi$
denote any character\mod q.

Since $\Zx{p^{\ord pq}}$ is a subgroup of $\Zqx$, we certainly have
$\ord p{E(q)} \ge \ord p{E(p^{\ord pq})} = \ord pq-1$ for all odd
primes $p$ dividing $q$, and also $\ord2{E(q)} \ge
\ord2{E(2^{\ord2q})} \ge \ord2q-2$. This implies that $\ord
p{E(q)/S(q)} \ge \max\{0,\ord pq-2\}$ for odd primes $p$ dividing $q$,
and $\ord2{E(q)/S(q)} \ge \max\{0,\ord2q-3\}$. In particular,
\begin{equation}
\ord p{\tilde{q_c} \cdot E(q)/S(q)} \ge \ord pq \quad \hbox{for all
$p\mid q$.} \label{qwert}
\end{equation}
If $p$ is a prime and $r$ a positive integer, we note that whenever
$m\equiv n\mod{p^r}$, we also have $m^p\equiv n^p\mod{p^{r+1}}$; one
can see this by using the binomial theorem to expand
$\{m+(n-m)\}^p-m^p$. Using this fact $\ord p{E(q)/S(q)}$ times for
each prime $p$ dividing $q$, and using the inequality (\ref{qwert}),
we can conclude that
\begin{equation*}
m\equiv n\mod{\tilde{q_c}} \implies m^{E(q)/S(q)}\equiv
n^{E(q)/S(q)} \mod q.
\end{equation*}
This tells us that for every $m\equiv n\mod{\tilde{q_c}}$, we have
\begin{equation*}
\chi^{E(q)/S(q)}(m) = \chi(m^{E(q)/S(q)}) = \chi(n^{E(q)/S(q)}) =
\chi^{E(q)/S(q)}(n),
\end{equation*}
which establishes the lemma.
\end{proof}

\begin{lemma}
For every $\chi\in G$, $N\ge1$, and $0<\eta<1$, we have
\begin{equation*}
\sum_{n\le N} \chi(n) \ll N \( N^{-1}q_c^{1/4+\eta} \)^\eta.
\end{equation*}
\label{Cqbd}
\end{lemma}

\begin{proof}
From Lemma \ref{induced}, $\chi$ is induced by a
character\mod{\tilde{q_c}}. However, $\tilde{q_c}$ is divisible by all
of the primes which divide $q$, and so $\chi$ is exactly equal to a
character\mod{\tilde{q_c}}. Write $\chi=\chi_1\chi_2$, where $\chi_1$
is a character\mod{q_c} and $\chi_2$ is a character\mod8. Then
\begin{equation*}
\begin{split}
\sum_{n\le N} \chi(n) &= \sum_{n\le N} \chi_1(n)\chi_2(n) \\
&= \sum_{i=1,3,5,7} \chi_2(i) \sum \begin{Sb}n\le N \\ n\equiv
i\mod8\end{Sb} \chi_1(n) \\
&\ll \max_{i=1,3,5,7} \abs{\sum_{m\le N/8} \chi_1(8m+i)}.
\end{split}
\end{equation*}
We choose $k$ so that $8k\equiv1\mod{q_c}$, and multiply the inner
sum through by $\chi_1(k)$; this doesn't change the size of the
expression, since $\abs{\chi_1(k)}=1$. Therefore
\begin{equation}
\sum_{n\le N} \chi(n) \ll \max_{j=k,3k,5k,7k} \abs{\sum_{j<n\le j+N/8}
\chi_1(n)}.
\label{maxjk}
\end{equation}
Now $\chi_1$ is a character\mod{q_c}, and $q_c$ is cubefree. By
Burgess' character sum estimates \cite[Theorem 2]{Bur:OCSaL-s2} for cubefree
moduli, for any integers $r\ge 2$ and $J$ and any $\ep>0$, we have
\begin{equation}
\begin{split}
\sum_{J<n\le J+N/8} \chi_1(n) &\ll (N/8)^{1-1/r}
q_c^{(r+1)/(4r^2)+\ep} \\
&\ll N \( N^{-1} q_c^{1/4+1/(4r)+\ep r} \)^{1/r}.
\end{split}
\label{rnoteta}
\end{equation}
If we stipulate that $r>\eta^{-1}$ and set $\ep=3\eta/(4r)$, the
bounds (\ref{maxjk}) and (\ref{rnoteta}) establish the lemma.
\end{proof}

\section{Sieve results}\label{sievesec}\noindent
\noindent In this section we cite the sieve results needed in the
later arguments. Most of the notation is standard for the theory of
sieves: let $\A$ be a set of positive integers, $X>1$ a real number,
and $\rho(d)$ a multiplicative function, and define the ``remainder
terms'' $R_d$ (which are intended to be small, at least on average
over $d$) by
\begin{equation*}
R_d = \sum \begin{Sb}a\in\A \\ d\mid a\end{Sb} 1 - \frac{\rho(d)}dX.
\label{Rddef}
\end{equation*}
Let $L$ and $y$ be positive numbers, and suppose that the
following conditions hold:
\begin{equation}
1 \ll 1-\frac{\rho(p)}p \le 1 \quad\text{for all primes }p;
\label{basicrhocond}
\end{equation}
\begin{equation}
-L < \sum_{w\le p<z} \frac{\rho(p)\log p}p - \log\frac zw < O(1)
\quad\text{for all }2\le w<z;
\label{advrhocond}
\end{equation}
\begin{equation}
\sum_{d\le y} \mu^2(d)3^{\omega(d)}\abs{R_d} \ll X/(\log X)^2.
\label{errorcond}
\end{equation}
Define $P(z)=\prod_{p<z}p$. Then we have
\begin{equation}
\sum \begin{Sb}a\in\A \\ (a,P(\sqrt y))=1\end{Sb} 1 \ll
X\prod_{p<\sqrt y}\( 1-\frac{\rho(p)}p \);
\label{upperlinear}
\end{equation}
also, for all $\ep>0$ and all $z\ge2$ such that $z^{2+\ep}\ll y$, we
have
\begin{equation}
\sum \begin{Sb}a\in\A \\ (a,P(z))=1\end{Sb} 1 \ge 
X\prod_{p<z}\( 1-\frac{\rho(p)}p \) \( \Delta(\ep) + O\(
{L\over(\log y)^{1/14}} \) \),
\label{lowerlinear}
\end{equation}
This follows, for example, from Theorems 4.1 and 8.4 of
Halberstam-Richert \cite{HalRic:SM}, where in Theorem 4.1, we take
$y=z^2$ in their notation, and in Theorem 8.4, we take $y=X^\alpha$ in
their notation and subsume the quantity $f(\cdot)$ into the constant
$\Delta(\ep)$.

We now describe the results of Greaves
on the weighted linear sieve that we will employ. Let $\A$, $X$,
$\rho$, $R_d$, $L$, and $y$ be as before (so that conditions
(\ref{basicrhocond}) through (\ref{errorcond}) are satisfied), and let
$g$ be a positive number such that
\begin{equation}
1\le a\le y^g \quad\text{for all }a\in\A. \label{ygcond}
\end{equation}
Let $r\ge2$ be an integer and $U$, $V$ be real numbers satisfying
\begin{equation}
V<1/4,\quad 1/2<U<1,\quad V+rU\ge g, \label{UVcond}
\end{equation}
and define $m=\max\{V,(1-U)/2\}$. Let $W$ be a nondecreasing ``weight
function'', defined on the positive real numbers, satisfying
\begin{equation}
0\le W(t)\le \begin{cases}
U-V &\rmif U<t, \\
t-V &\rmif 1/3\le t\le U, \\
t-m &\rmif m<t\le1/3 \\
\end{cases}
\label{Wcond}
\end{equation}
and also
\begin{equation*}
W(t)\le9(U-1/3)t^2\quad\text{for all }t>0. \label{Wcond2}
\end{equation*}
Define $\omega(z,n)$ to be the number of prime factors of $n$, where a
multiple prime factor $p$ is counted multiply only if $p\ge z$, so
that $\omega(\infty,n)=\omega(n)$ and $\omega(2,n)=\Omega(n)$. Also
define $p(n)$ to be the least prime factor of $n$. Then
\begin{equation}
\sum \begin{Sb}a\in\A \\ \omega(y^U,a)\le r\end{Sb} W\( \log
p(a)\over\log y \) \ge 2e^\gamma X\prod_{p<y}\( 1-\frac{\rho(p)}p \) \(
{\cal M}(W) + O\( \( L\over\log y\)^{1/5} \) \),
\label{Gresult}
\end{equation}
where $\gamma$ is Euler's constant and ${\cal M}(W)$ is a constant
depending on $W$. This is essentially Theorem 1 of Greaves
\cite{Gre:TWLSaSetc.}, although we have already dealt with the
remainder term through the condition (\ref{errorcond}).

Furthermore, there exist $U$, $V$ satisfying condition (\ref{UVcond})
and a weight function $W$ satisfying conditions (\ref{Wcond}) and
(\ref{Wcond2}) such that ${\cal M}(W)$ is positive, as long as $g\le
r-\delta_r$, where the values in Theorem \ref{mainthm} are permissible
values for $\delta_r$. (The permissible value when $r\ge5$ comes from
Greaves' earlier work \cite{Gre:AWSoBT}.)

\section{Proof of Theorem \ref{mainthm}}\noindent
\noindent The following lemma provides an upper bound for a remainder
term sum we will encounter in our applications of the sieve results in
Section \ref{sievesec}. We recall that $c_0$ is a shorthand for
$c_{\chi_0}$.

\begin{lemma}
Given a real number $x>1$ and coprime integers $q\ge2$ and $d$
satisfying $1\le d<x$, let the quantity $R_d$ be defined by
\begin{equation}
R_d = \sum \begin{Sb}n<x \\ d|n\end{Sb} \gamma(n) - \frac{c_0\phi(q)}q
\frac xd,
\label{Rdform}
\end{equation}
Then for any $\ep$, $\eta>0$, we have
\begin{equation*}
R_d \ll \frac{c_0\phi(q)}q\frac xdq_c^{2\ep}
\( \frac dxq_c^{1/4+\eta} \)^\eta.
\end{equation*}
\label{remainlem}
\end{lemma}

\begin{proof}
Using Lemma \ref{charlincomb}, we can write
\begin{equation}
\begin{split}
R_d &= \sum\begin{Sb}n<x \\ d\mid n\end{Sb} \sum_{\chi\in G}
c_\chi\chi(n) - \frac{c_0\phi(q)}q \frac xd \\
&= \sum_{\chi\in G} c_\chi\chi(d) \sum_{m<x/d} \chi(m) -
\frac{c_0\phi(q)}q \frac xd.
\end{split}
\label{absAd}
\end{equation}
For the term corresponding to the principal character, we note that
for any $T>1$, we have
\begin{equation*}
\begin{split}
\sum_{n<T} \chi_0(n) &= \sum \begin{Sb}n<T \\ (n,q)=1\end{Sb} 1 =
\sum_{n<T} \sum \begin{Sb}f\mid n \\ f\mid q\end{Sb} \mu(f) =
\sum_{f\mid q} \mu(f) \sum_{m<T/f} 1 \\
&= T \sum_{f\mid q} \frac{\mu(f)}f + O\( \sum_{f\mid q} \abs{\mu(f)}
\) \\
&= T\frac{\phi(q)}q + O\( 2^{\omega(q)} \).
\end{split}
\end{equation*}
Thus
\begin{equation}
c_0\chi_0(d)\sum_{m<x/d}\chi(m) = \frac{c_0\phi(q)}q\frac xd 
+ O\( c_0\cdot2^{\omega(q)} \),
\label{principal}
\end{equation}
the first term of which will cancel the last term of equation
(\ref{absAd}). For the other terms in the sum over $\chi$ in
(\ref{absAd}), we apply Lemma \ref{Cqbd} to the inner sums to see that
\begin{equation}
\sum_{m<x/d} \chi(m) \ll \frac xd \( \frac dxq_c^{1/4+\eta}
\)^\eta
\label{otherchis}
\end{equation}
for any $\ep$, $\eta>0$. We use equations (\ref{principal}) and
(\ref{otherchis}) in equation (\ref{absAd}) to get
\begin{equation*}
\begin{split}
R_d &\ll 2^{\omega(q)}c_0 + \sum
\begin{Sb}\chi\in G \\ \chi\ne\chi_0\end{Sb} \abs{c_\chi}
\frac xd \( \frac dxq_c^{1/4+\eta} \)^\eta \\
&\ll \frac{c_0\phi(q)}q \frac xdq_c^{2\ep}
\( \frac dxq_c^{1/4+\eta} \)^\eta,
\end{split}
\label{absAd2}
\end{equation*}
since $\sum_{\chi\in G}\abs{c_\chi}=2^{\omega(\phi(q))}c_0$ by Lemma
\ref{cchisum}, and both $q/\phi(q)$ and $2^{\omega(\phi(q))}$ are $\ll
q_c^\ep$. This establishes the lemma.
\end{proof}

\begin{lemma}
Let $q$, $x$, and $\eta$ be as in Lemma \ref{remainlem}, and define
\begin{equation}
y=x^{1-\eta}/q_c^{1/4+3\eta}. \label{ydef}
\end{equation}
Then
\begin{equation*}
\sum_{d\le y} \mu^2(d)3^{\omega(d)}R_d \ll \frac{c_0\phi(q)}q
x^{1-\eta^2/2}.
\end{equation*}
In particular, condition (\ref{errorcond}) is satisfied.
\label{remaincor}
\end{lemma}

\begin{proof}
Let $\ep=\eta^2$. For any $n<x$, we have $3^{\omega(n)}\ll x^{\ep/2}$,
and so
\begin{equation}
\sum_{d\le y} \mu^2(d)3^{\omega(d)}R_d \ll x^{\ep/2} \sum_{d\le
y} \abs{R_d}. \label{triangle}
\end{equation}
When $(d,q)>1$, we have $R_d=0$; and so we may use equation
(\ref{triangle}), Lemma \ref{remainlem}, and the definition
(\ref{ydef}) of $y$ to see that
\begin{equation*}
\begin{split}
\sum_{d\le y} \mu^2(d)3^{\omega(d)}R_d &\ll x^{\ep/2}
\sum_{d\le y} \frac{c_0\phi(q)}q\frac xdq_c^{2\ep}
\( \frac dxq_c^{1/4+\eta} \)^\eta \\
&\ll x^{\ep/2} \frac{c_0\phi(q)}q xq_c^{2\ep} \(
\frac yxq_c^{1/4+\eta} \)^\eta \\
&= x^{\ep/2} \frac{c_0\phi(q)}q xq_c^{2\ep}
(x^{-\eta}q_c^{-2\eta})^\eta = \frac{c_0\phi(q)}q x^{1-\eta^2/2},
\end{split}
\end{equation*}
which establishes the lemma.
\end{proof}

We now use the weighted linear sieve to deduce the following
quantitative version of Theorem \ref{mainthm}.

\begin{theorem}
Let $q$, $r\ge2$ be integers, $0<\eta<1/12$ a real number, and $x$ a
real number satisfying
\begin{equation}
x \ge q_c^{1/4+1/(4(r-1-\delta_r))+15\eta}, \label{xbigenough}
\end{equation}
where the $\delta_r$ are given in Theorem \ref{mainthm}. Then for some
positive constant $C=C(r,\eta)$, we have
\begin{equation*}
\sum \begin{Sb}n<x \\ n=P_r\end{Sb} \gamma(n) \ge \frac{c_0\phi(q)}q
\frac x{\log x} \( C + O\( \log\log 3x\over\log x \)^{1/5} \).
\end{equation*}
\label{mainqthm}
\end{theorem}


\begin{proof}
We apply the lower bound (\ref{Gresult}), with $\A$ being the set of
all \lr{s}\mod q less than $x$, and with the various parameters taken
as follows:
\begin{equation}
\begin{array}{c}
\displaystyle X = \frac{c_0\phi(q)}qx ; \qquad g = r-\delta_r; \qquad
\rho(p) = \begin{cases}
1 &\rmif p\dnd q, \\ 0 &\rmif p\mid q;
\end{cases}
\vphantom{\Bigg(}\\
\displaystyle y = x^{1-\eta}/q_c^{1/4+3\eta}; \qquad L \ll
\log\log3q_c.
\end{array}
\label{defs}
\end{equation}
We note that for $d$ coprime to $q$, the remainder term $R_d$ takes
the form given by equation (\ref{Rdform}), while for $d$ not coprime
to $q$, we have $R_d=0$. It is straightforward to verify the
conditions (\ref{basicrhocond}), (\ref{advrhocond}), and
(\ref{ygcond}) for this choice of parameters, the requirement
(\ref{xbigenough}) being crucial to the validity of
(\ref{ygcond}). Moreover, we see from Lemma \ref{remaincor} that
condition (\ref{errorcond}) holds as well. Therefore, equation
(\ref{Gresult}) gives us
\begin{equation}
\sum \begin{Sb}n\le x \\ \omega(y^U,n)\le r\end{Sb} \gamma(n) W\( \log
p(n)\over\log y \) \ge \frac{c_0\phi(q)}q \frac x{\log x} \( {\cal
M}(W) + O\( \( \log\log3x\over\log x\)^{1/5} \) \).
\label{application}
\end{equation}
where we have used Mertens' formula for the product over primes less
than $y$ and the fact that, from the definition (\ref{defs}) of $y$
and the restriction (\ref{xbigenough}) on $x$, we have
\begin{equation*}
\log x\ge\log y\gg\log x \quad\text{and}\quad \log\log q_c\ll\log\log x.
\end{equation*}
The terms counted by the sum in (\ref{application}) are not
necessarily $P_r$s; we now assure that the contribution of those
integers that are not $P_r$s is negligible. First we estimate the
order of magnitude of the factor $c_0\phi(q)/q$. By the definition of
$c_0$, we see that $c_0
\ge \phi(\phi(q))/\phi(q)$. Since $n/\phi(n)\ll\log\log n$ for
any positive integer $n$, the restriction (\ref{xbigenough}) implies
that
\begin{equation*}
q/(c_0\phi(q)) \ll (\log\log3\phi(q))(\log\log3q) \ll (\log\log3x)^2.
\end{equation*}
Notice that, by the condition (\ref{Wcond}) on $W$ and the fact that
$W$ is nondecreasing, any integer counted with a positive weight by
the sum in (\ref{application}) has no prime factors less than $y^m$.
Since $W(t)\le W(1)$ for all $t>0$, we may define $C={\cal M}(W)/W(1)$
and divide both sides of equation (\ref{application}) by $W(1)$ to see
that
\begin{equation*}
\sum \begin{Sb}n\le x \\ \omega(y^U,n)\le r \\ p(n)>y^m\end{Sb}
\gamma(n) \ge \frac{c_0\phi(q)}q \frac x{\log x} \( C + O\( \(
\log\log3x\over\log x \)^{1/5} \) \).
\end{equation*}
Now the number of squarefree integers less than $x$ whose smallest
prime factor exceeds $y^m$ is at most
\begin{equation*}
\sum_{p>y^m} \frac x{p^2} \ll \frac x{y^m}.
\end{equation*}
Thus
\begin{equation}
\sum \begin{Sb}n\le x \\ \omega(y^U,n)\le r \\ p(n)>y^m\end{Sb}
\mu^2(n)\gamma(n) \ge \frac{c_0\phi(q)}q \frac x{\log x} \( C + O\( \(
\log\log3x\over\log x \)^{1/5} + {(\log x)(\log\log3x)^2\over y^m} \) \),
\label{finally}
\end{equation}
and again the second error term is negligible, since the definition
(\ref{defs}) of $y$ and the restriction (\ref{xbigenough}) on $x$
insure (for $\eta<\frac13$, for instance) that $y\ge x^\eta$. But if $n$
is squarefree, then $\omega(z,n)=\Omega(n)$ for any real $z$. Thus the
only integers counted by the sum in (\ref{finally}) are $P_r$s, which
establishes the theorem.
\end{proof}

\section{Proof of Theorem \ref{primrootthm}}\label{rificsec}\noindent
Theorem \ref{primrootthm} is almost a corollary of Theorem
\ref{mainqthm}, except that we must argue that some of the small
almost-prime primitive roots\mod p are also primitive roots\mod{p^2}.
The following lemma was established by Kruswijk \cite{Kru:OtCetc}; we
provide a proof for the sake of completeness.
%

\begin{lemma}
Let $p$ be a prime, and for every real $x>1$, define the set
\begin{equation*}
\B(x) = \{b\le x: \text{\rm$b$ is a \pth\ power\mod{p^2}}\}.
\end{equation*}
Let $B(x)=\abs{\B(x)}$. Then, uniformly for all positive integers $m$,
we have
\begin{equation*}
B(p^{1/m}) \le p^{1/(2m)} \exp\( O\( \log p\over\log\log p \) \).
\end{equation*}
\label{deviouslem}
\end{lemma}

\begin{proof}
Fix an integer $m\ge1$ and consider the set
\begin{equation*}
\C = \{b_1\cdots b_{2m}: b_i\in\B(p^{1/m}),\,1\le i\le2m\}.
\end{equation*}
On one hand, every element of $\C$ is a \pth\ power\mod{p^2} and is at
most $p^2$ in size; thus $\abs\C\le p$, the total number of \pth\
powers\mod{p^2} up to $p^2$. On the other hand, the total number of
products of $2m$ elements of $\B(x)$ is of course
$B(x)^{2m}$. Moreover, the number of ways to write an integer $n$ as a
product of $2m$ elements of $\B(x)$ is bounded by the number of ways
to decompose $n$ generally as a product of $2m$ factors; and this in
turn is at most
\begin{equation*}
(2m)^{\Omega(n)},
\end{equation*}
as we can assign each of the $\Omega(n)$ prime divisors of $n$ to any
of the $2m$ factors in the decomposition, and this covers all
cases. When we take $x=p^{1/m}$, the fact that $n\le p^2$ for all
$n\in\C$ together with a standard upper bound for $\Omega(n)$ gives
us
\begin{equation*}
\Omega(n) \ll {\log p^2\over\log\log p^2}.
\end{equation*}
Thus we can say that
\begin{equation*}
p \ge \abs\C \ge {B(p^{1/m})^{2m} \over (2m)^{O(\log p/\log\log p)}},
\end{equation*}
or
\begin{equation*}
B(p^{1/m}) \le p^{1/(2m)} \exp\( O\( {\log m\over m} {\log
p\over\log\log p} \) \).
\end{equation*}
Since $\log m/m$ is uniformly bounded, this establishes the lemma.
\end{proof}

We can now establish Theorem \ref{primrootthm}. First, when $q=p$ is
an odd prime, Theorem \ref{primrootthm} is an immediate corollary of
Theorem \ref{mainqthm}, using the values of $\delta_r$ given in
Theorem \ref{mainthm}. Now the only way a primitive root\mod p can
fail to be a primitive root\mod{p^2} is if it is a \pth\
power\mod{p^2}. Applying Lemma \ref{deviouslem} with various values of
$m$, we see that
\begin{equation*}
\begin{split}
B(p^{1/2+1/873}) &\le B(p) \ll p^{1/2+\ep}; \\
B(p^{1/3+1/334}) \le B(p^{3/8+1/207}) &\le B(p^{1/2}) \ll p^{1/4+\ep};
\\
B(p^{1/4+O(1/r)}) &\le B(p^{1/3}) \ll p^{1/6+\ep}.
\end{split}
\end{equation*}
In all cases, the number of \pth\ powers\mod{p^2} is of a lower order
of magnitude than the number of primitive roots\mod p. This
establishes Theorem \ref{primrootthm} for $q=p^2$. Finally, it is
well-known that any primitive root\mod{p^2} is also a primitive
root\mod{p^r} for every $r\ge3$, and that any odd primitive
root\mod{p^r} is also a primitive root\mod{2p^r} for every
$r\ge1$. Since the primitive roots counted by Theorem \ref{mainqthm}
are odd, this establishes Theorem \ref{primrootthm} in its entirety.

\section{Siegel zeros and prime primitive roots: Proof of Theorem
\ref{Siegelthm}}\noindent
In this section we only consider moduli $q$ that admit primitive
roots, i.e., $q$ is an odd prime power or twice an odd prime power (we
ignore $q=2$ and $q=4$). Every primitive root\mod q is certainly a
quadratic nonresidue. Thus those primitive roots with an even number
of prime factors (counted with multiplicity) must have a prime factor
that is a quadratic residue.

On the other hand, suppose that for such a modulus $q$, the
$L$-function of the unique nonprincipal quadratic character $\chi_1$
has a Siegel zero. We know, by the prime number theorem for arithmetic
progessions, that this makes small primes that are quadratic residues
very rare, and so small primitive roots with an even number of prime
factors are correspondingly rare. Thus if we can show the existence of
$P_2$ primitive roots of a certain type, then assuming the existence
of a Siegel zero, we might expect to be able to argue that most of
them must in fact be primes.

To do so, we use the lower bound linear sieve to produce $P_2$
primitive roots, and then a simple upper bound sieve to show that the
contribution from primitive roots divisible by a quadratic residue is
small. First we need to consider the remainder sum we will encounter
while employing the upper bound sieve.

\begin{lemma}
Let $q$ be an odd prime and $x>1$ a real number, and let $\chi_1$ be
the nonprincipal quadratic character\mod q. Define the quantity
\begin{equation*}
H = \sum\begin{Sb}x^{1/3}<p<x^{2/3} \\ \chi_1(p)=1\end{Sb} \frac1p.
\label{Hdef}
\end{equation*}
Given an integer $d$, coprime to $q$ and satisfying $1\le d<x^{1/3}$,
let the quantity $R_d$ be defined by
\begin{equation}
R_d = \sum\begin{Sb}x^{1/3}<p<x^{2/3} \\ \chi_1(p)=1\end{Sb} 
\sum\begin{Sb}n<x/p \\ d\mid n\end{Sb} \gamma(pn) -
\frac{c_0\phi(q)}q\frac xdH,
\label{RddefH}
\end{equation}
Then for any $\ep$, $\eta>0$, we have
\begin{equation*}
R_d \ll \frac{\phi(q-1)}q\frac xdHq^{2\ep}
\( \frac d{x^{1/3}}q^{1/4+\eta} \)^\eta.
\end{equation*}
\label{remainlemH}
\end{lemma}

\begin{proof}
Since the argument parallels the proof of Lemma \ref{remainlem}, we
provide only an outline of the proof. We have
\begin{equation*}
R_d = \sum_{\chi\in G} c_\chi\chi(d) \sum\begin{Sb}x^{1/3}<p<x^{2/3}
\\ \chi_1(p)=1\end{Sb} \chi(p) \sum_{m<x/pd} \chi(m) -
\frac{\phi(q-1)}q\frac xdH.
\end{equation*}
For the principal character we have
\begin{equation*}
\begin{split}
c_0\chi_0(d)\sum\begin{Sb}x^{1/3}<p<x^{2/3} \\ \chi_1(p)=1\end{Sb}
\chi_0(p) \sum_{m<x/pd} \chi_0(m) - \frac{\phi(q-1)}q\frac xdH &\ll
{\phi(q-1)\over q-1} \sum\begin{Sb}x^{1/3}<p<x^{2/3} \\
\chi_1(p)=1\end{Sb} 2^{\omega(q-1)} \\
&\ll {\phi(q-1)\over q-1} 2^{\omega(q-1)}x^{2/3}H,
\end{split}
\end{equation*}
while for $\chi\ne\chi_0$ we have
\begin{equation*}
\begin{split}
c_\chi\chi(d)\sum\begin{Sb}x^{1/3}<p<x^{2/3} \\ \chi_1(p)=1\end{Sb}
\chi(p) \sum_{m<x/pd} \chi(m) &\ll
\abs{c_\chi}\sum\begin{Sb}x^{1/3}<p<x^{2/3} \\ \chi_1(p)=1\end{Sb}
\frac x{dp} \( \frac{dp}xq^{1/4+\eta} \)^\eta \\
&\ll \abs{c_\chi}\frac xdH \( \frac d{x^{1/3}}q^{1/4+\eta} \)^\eta.
\end{split}
\end{equation*}
Thus
\begin{equation*}
R_d \ll \frac{\phi(q-1)}q\frac xdHq^{2\ep} \( \frac
d{x^{1/3}}q^{1/4+\eta} \)^\eta,
\end{equation*}
which establishes the lemma.
\end{proof}

\begin{lemma}
Let $q$, $x$, and $\eta$ be as in Lemma \ref{remainlemH}, and define
\begin{equation*}
y=x^{1/3-\eta}/q^{1/4+3\eta}. \label{ydefH}
\end{equation*}
Then
\begin{equation*}
\sum_{d\le y} \mu^2(d)3^{\omega(d)}R_d \ll \frac{\phi(q-1)}q
x^{1-\eta^2/2}H.
\end{equation*}
In particular, condition (\ref{errorcond}) is satisfied.
\label{remaincorH}
\end{lemma}

\begin{proof}
Let $\ep=\eta^2$. For any $n<x$, we have $3^{\omega(n)}\ll x^{\ep/2}$,
and so
\begin{equation*}
\sum_{d\le y} \mu^2(d)3^{\omega(d)}R_d \ll x^{\ep/2} \sum_{d\le
y} \abs{R_d}.
\end{equation*}
We may now apply Lemma \ref{remainlemH} to see that
\begin{equation*}
\begin{split}
\sum_{d\le y} \mu^2(d)3^{\omega(d)}\abs{R_d} &\ll x^{\ep/2} \sum_{d\le
y} \frac{\phi(q-1)}q\frac xdHq^{2\ep} \( \frac d{x^{1/3}}q^{1/4+\eta}
\)^\eta \\
&\ll x^{\ep/2} \frac{\phi(q-1)}q xHq^{2\ep} \( \frac
y{x^{1/3}}q^{1/4+\eta} \)^\eta \\
&\ll x^{\ep/2} \frac{\phi(q-1)}q xHq^{2\ep} (x^{-\eta}q^{-2\eta})^\eta
= \frac{\phi(q-1)}q x^{1-\eta^2/2}H,
\end{split}
\end{equation*}
which establishes the lemma.
\end{proof}

We may now establish the following quantitative version of Theorem
\ref{Siegelthm}.

\begin{theorem}
Let $q$ be an odd prime power or twice an odd prime power, and let
$\chi_1$ be the nonprincipal quadratic character\mod q. Suppose
that $L(s,\chi_1)$ has a real zero $\beta$ of the form
\begin{equation*}
\beta = 1-\frac1{\alpha\log q}
\end{equation*}
with $\alpha\ge3$. Then for any real number $0<\eta<1/52$, and any real
number $x$ satisfying $q^{3/4+13\eta}\le x\le q^{500}$, we have
\begin{equation*}
\sum_{p<x} \gamma(p) \ge \frac{\phi(\phi(q))}q \frac x{\log x} \(
\Delta(\eta) + O\( (\log x)^{-1/14} + (\log\alpha)^{-1/2} \) \),
\end{equation*}
where $\Delta(\eta)$ is a positive constant depending on $\eta$.
\label{quantSiegel}
\end{theorem}

\begin{proof}
First we assume that $q$ is an odd prime. We let $\A$ be the set of
all primitive roots\mod q less than $x$, and we choose
\begin{equation}
\begin{array}{c}
\displaystyle X = \frac{\phi(q-1)}qx; \qquad \rho(d) = 1\text{ for all
}d<q;
\vphantom{\Bigg(}\\
\displaystyle y = x^{1-\eta}/q^{1/4+3\eta}; \qquad z = x^{1/3}; \qquad
L \ll1.
\end{array}
\label{defsagain}
\end{equation}
Again, conditions (\ref{basicrhocond}) and (\ref{advrhocond}) are easy
to verify, and the special case of Lemma \ref{remaincor} where $q$ is
an odd prime establishes that condition (\ref{errorcond}) holds. In
addition, the definition (\ref{defsagain}) of $y$, together with the
restrictions $x\ge q^{3/4+13\eta}$ and $\eta<\frac1{52}$, insure that
$z^{2+\eta}\le y$. Therefore we may apply the linear sieve
(\ref{lowerlinear}) to obtain
\begin{equation}
\sum \begin{Sb}a\in\A \\ (a,P(x^{1/3}))=1\end{Sb} 1 \ge 
\frac{\phi(q-1)}q \frac x{2e^\gamma\log x^{1/3}} \( \Delta(\eta) + O\(
(\log x)^{-1/14} \) \),
\label{unweighted}
\end{equation}
where we have again used Mertens' formula for the product over
primes. We notice that all the integers counted by this sum are
$P_2$s, and that the lower bound is all the more true if we replace
the denominator $2e^\gamma\log x^{1/3}$ by $\log x$.

We now show that the contribution to the sum (\ref{unweighted}) from
products of two primes is negligible. We write this contribution as
\begin{equation*}
T = \sum\begin{Sb}p_1p_2<x \\ p_1,p_2>x^{1/3}\end{Sb}
\gamma(p_1p_2).
\label{SminusS}
\end{equation*}
When $p_1p_2$ is a primitive root, which must be the case for $p_1p_2$
to contribute to this sum, exactly one of $p_1$ or $p_2$ is a
quadratic residue\mod q. Thus this sum becomes
\begin{equation*}
T = \sum\begin{Sb}x^{1/3}<p_1<x^{2/3} \\ \chi_1(p_1)=1\end{Sb}
\sum\begin{Sb}x^{1/3}<p_2<x/p_1 \\ \chi_1(p_2)=-1\end{Sb}
\gamma(p_1p_2).
\end{equation*}
Relaxing the restrictions on the integer $p_2$ will allow us to
estimate this sum more easily. Let $z$ be a parameter satisfying
$2\le z\le x^{1/3}$. Instead of summing over primes in the range
$x^{1/3}<p<x/p_1$ which are quadratic nonresidues, we will simply sum
over all integers $n<x/p_1$ whose prime divisors all exceed $z$. Thus
we have
\begin{equation}
T \le T(z) = \sum\begin{Sb}x^{1/3}<p<x^{2/3} \\ \chi_1(p)=1\end{Sb}
\sum\begin{Sb}n<x/p \\ (n,P(z))=1\end{Sb} \gamma(pn).
\label{Tdef}
\end{equation}
We can now apply the upper bound sieve (\ref{upperlinear}) with $\A$
being the set of all integers $n<x^{2/3}$, counted with multiplicity
equal to the number of primes $p$, meeting the criteria of the outer
sum in (\ref{Tdef}), such that $pn$ is a primitive root\mod q. We
take the various sieve parameters as follows:
\begin{equation*}
\begin{array}{c}
\displaystyle X = \frac{c_0\phi(q)}qxH; \qquad \rho(d) = 1\text{ for
all }d<q;
\vphantom{\Bigg(}\\
\displaystyle y = x^{1/3-\eta}/q^{1/4+3\eta}; \qquad L\ll1.
\end{array}
\end{equation*}
Now the remainder term $R_d$ takes the form given by equation
(\ref{RddefH}), and it is again straightforward to verify the
conditions (\ref{basicrhocond}) and (\ref{advrhocond}) and to note
that Lemma \ref{remaincorH} establishes condition (\ref{errorcond}) as
well. Therefore, equation (\ref{upperlinear}) gives us
\begin{equation*}
T \le T(\sqrt y) \ll \frac{\phi(q-1)}q\frac x{\log x}H.
\end{equation*}
which, together with the lower bound (\ref{unweighted}), yields
\begin{equation}
\sum \begin{Sb}p\in\A \\ p>x^{1/3}\end{Sb} 1 \ge 
\frac{\phi(q-1)}q \frac x{\log x} \( \Delta(\eta) + O\(
(\log x)^{-1/14} + H \) \),
\label{stillH}
\end{equation}
Heath-Brown has shown \cite[Lemma 3]{Hea:PTaSZ} that for $x\le
q^{500}$, we have
\begin{equation*}
\sum \begin{Sb}p<x \\ \chi_1(p)=1 \end{Sb} \frac{\log p}p \ll (\log
q)(\log\alpha)^{-1/2},
\end{equation*}
from which it follows that $H \ll (\log\alpha)^{-1/2}$, since $\log
x\gg\log q$. Using this fact in the lower bound (\ref{stillH})
establishes the theorem when $q$ is an odd prime. The arguments used
in Section \ref{rificsec} to justify Theorem \ref{primrootthm} for
composite moduli apply here as well to establish the theorem in its
full generality.
\end{proof}

\bibliography{primroots}
\bibliographystyle{amsplain}
\end{document}